\documentclass[journal]{IEEEtran}
\IEEEoverridecommandlockouts
\usepackage{cite}
\usepackage{amsmath,amssymb,amsfonts}
\usepackage{algorithmic}
\usepackage{graphicx}
\usepackage{textcomp}
\usepackage{xcolor}
\usepackage{bm}
\usepackage{lineno,hyperref,float,amsfonts,amsmath,epsfig,graphics,booktabs,multirow,cite,stfloats}
\usepackage{amsthm,mathrsfs}
\usepackage{url}
\usepackage{colortbl}
\usepackage{subfigure}
\usepackage{verbatim}
\usepackage{comment}

\usepackage{algorithm}

\theoremstyle{definition}

\theoremstyle{plain}

\theoremstyle{remark}

\def\BibTeX{{\rm B\kern-.05em{\sc i\kern-.025em b}\kern-.08em
    T\kern-.1667em\lower.7ex\hbox{E}\kern-.125emX}}
\begin{document}

\title{Fine-tuning for Data-enabled Predictive Control of Noisy Systems by Reinforcement Learning}
\author{Jinbao Wang, Shiliang Zhang, Jun Liu, Xuehui Ma, Haolin Liu
\thanks{J. Wang, J. Liu, X. Ma, and H. Liu are with Xi’an University of Technology, Xi'an, China (email: jbwang@stu.xaut.edu.cn, liujun0310@xaut.edu.cn, xuehui.yx@gmail.com, haolinliu@xaut.edu.cn, S. Zhang is with University of Oslo, Norway (email: shilianz@uio.no)}
}

\maketitle

\begin{abstract}

%% Text of abstract

Data-enabled predictive control (DeePC) leverages system measurements in characterizing system dynamics for optimal control. The performance of DeePC relies on optimizing its hyperparameters, especially in noisy systems where the optimal hyperparameters adapt over time. Existing hyperparameter tuning approaches for DeePC are more than often computationally inefficient or overly conservative. This paper proposes an adaptive DeePC where we guide its hyperparameters adaption through reinforcement learning. We start with establishing the relationship between the system I/O behavior and DeePC hyperparameters. Then we formulate the hyperparameter tuning as a sequential decision-making problem, and we address the decision-making through reinforcement learning. We implement offline training to gain a reinforcement learning model, and we integrate the trained model with DeePC to adjust its hyperparameters adaptively in real time. We conduct numerical simulations with diverse noisy conditions, and the results demonstrate the identification of near-optimal hyperparameters and the robustness of the proposed approach against noises in the control

\end{abstract}


\section{Introduction}
The control of uncertain systems with noise has been a challenging issue that triggered extensive studies over recent decades~\cite{DataInformativity,yu2024adaptive, BacksteppingControl,10591250,MPC,10551447,ma2022adaptive,ma2022adaptivenew,ma2022active,ma2020dual}. Those studies can be categorized into model-based and model-free control strategies. Model-based approaches, \textit{e.g.}, stochastic MPC~\cite{zhang2017data,yuan2016nonlinear,zhang2016nonlinear,zhang2016energy,zhang2015predictive}, use noisy data in identifying unknown model parameters and subsequently deriving the control law~\cite{zhang2021distributed,eslami2024control}. However, gaining the accurate model from contaminated measurements is nontrivial especially for large-scale complex systems~\cite{othermatrix} and the conditions for parameter identification are more than often case-specific. In this paper, we look into the model-free approaches, which directly map noisy data to optimal control signals without explicit model parameter recognition~\cite{Data-Driven, Data-Informativity}. Particularly, we examine data-enabled predictive control (DeePC), which is model free approach widely applied in uncertain and complex real systems like vehicle applications~\cite{liu2025adaptive,10901967,bai2025exploring,yang2024hardware,bai2024long}, energy applications~\cite{zhang2024impact,zhang2022extended,zhang2021energy}, and power electronic systems~\cite{Willems,Approximate,Ivan2008,Ivan2017,MARKOVSKY202142,markovsky2022data}.

DeePC is a model-free control rooted in Willems' fundamental lemma. It adopts regularized optimal control to ensure system performance and safety under uncertain conditions, such as data corruption, noises~\cite{DeePCNonlinear,DeePCPowerSystem,DeePCMotor,SegmentedTrajectories,AlleviatingThermal}, and outliers~\cite{zhang2019outlier,8291826}. DeePC uses data to estimate system dynamics and optimize the control, and preset hyperparameters for the cost function, the optimization of which leads to optimal control. 

DeePC's performance relies much on its hyperparameter, particularly the regularization term $\lambda_g$ that balances different objectives in the cost function. Huang \textit{et al.} demonstrated that $\lambda_g$ plays a critical role in the robustness of DeePC~\cite{huang2023robust}. Mattsson \textit{et al.} revealed that maintaining a fixed $\lambda_g$ can lead to declined control performance as the data volume increases~\cite{mattsson2024equivalence}. O'Dwyer \textit{et al.} pointed out that real-world DeePC implementations demand a rigorous tuning methodology to guarantee validity of the control~\cite{o2022data}. Real time adjustment and optimization of $\lambda_g$ is essential to ensure optimal performance~\cite{PowerSystems}. However, though DeePC has been well investigated, a solid DeePC tuning for uncertain systems remains missing. Below we analyze relevant parameter adjustment methods for DeePC and check the merits and limitations.

Existing hyperparameter tuning approaches for DeePC can be classified into three categories. The first category entails empirical assignment of hyperparameters. Such a method ignores system dynamics and can compromise the control~\cite{fiedler2021relationship,li2024physics,dai2022cloud}. The second category conducts an exhaustive search for the optimal value of~\( \lambda_g \) \cite{PowerSystems,alpago2020extended,van2024subspace,huang2021decentralized} in DeePC tuning. Although this method means to find the optimal or near-optimal parameter, the gained parameters are dedicated to the system of interest and the associated conditions. As a result, even minor changes in the system or its surroundings might lead to degraded control. 

The third category concerns approaches for real-time estimation and adjustment of the DeePC hyperparameter. Among those approaches, Chiuso \textit{et al.} imposed a dynamic penalty on the regularization term~\cite{CHIUSO2025112070} in tuning DeePC. They tune the regularization term $\lambda_g$ in real time based on the deviation between predicted and actual system outputs. However, their approach heavily depends on prior knowledge about the system model. Huang \textit{et al.} reformulated the hyperparameter selection of \( \lambda_g \) as a constrained optimization problem~\cite{huang2023robust}. While theoretically sound, this approach entails recurrent optimization calculations that risk significant overhead. Lazar \textit{et al.} recast the regularized DeePC framework as a Tikhonov-regularized least squares formulation, where they employ Hanke-Raus rule for parameter tuning~\cite{lazar2022offset}. Though their approach obviates extensive manual hyperparameter tuning, the generalization of their approach remains uncertain across different system configurations~\cite{verheijen2023handbook}. Furthermore, the analyzed approaches lack an applicable strategy that can be adapted to varying and uncertain system dynamics. In addressing this issue, Cummins \textit{et al.} introduced a backpropagation-based optimization framework, wherein the DeePC hyperparameter is fine-tuned using an approximate system model~\cite{cummins2024deepc}. While this approach allows for gradient-based optimization that offers adaptiveness for DeePC hyperparameter tuning, it induces significant computations and relies on accurate system approximation, and such a strict condition may not necessarily hold in practice.

In this paper, we provide real-time fine-tuning for DeePC in the control of uncertain systems with noises and varying system parameters. We base our design on no prior knowledge about the system, and we adjust the regularization term $\lambda_g$ according to the estimated system dynamics that better fit the control. To this end, we analyze the hyperparameter's impact on the system input-output (I/O) dynamics, and explore the relationship between environmental noise and the hyperparameter. Upon the gained relationship, we develop a reinforcement learning approach that adjust DeePC hyperparameter that coordinates the uncertainties and serves the control best. While we dedicate our fine-tuning approach to DeePC, it is possible to extend our approach to control strategies where a regularization term is necessary in balancing objectives in the control and demands real time adjustment to fit the system dynamics. We summarize the contributions of this study in the following.

\begin{itemize}
	\item[1] We systematic analyze the relationship between DeePC hyperparameter \( \lambda_g \) and the system's input-output (I/O) behavior in revealing the interactions between the DeePC controller and the varying system dynamics.
	\item[2] We introduce a reinforcement-learning based fine-tuning for DeePC, which adaptively adjusts the DeePC hyperparameter without model parameter identification or any prior knowledge. The developed fine-tuning can also be leveraged in other control systems where real time hyperparameter adjustment is necessary.
    \item[3] The fine-tuning better facilitates DeePC in the control of uncertain systems with varying dynamics and noises, without considerable increase in the computation that risks overhead. We use offline training for the reinforcement learning, such that we integrate the trained model to DeePC to make real-time decisions on hyperparameter tuning during the control.
\end{itemize}

The remainder of this paper is as follows. Section 2 formulates the hyperparameter tuning problem in the DeePC framework. Section 3 analyzes the relationship between the hyperparameter \( \lambda_g \) and the system’s I/O characteristics and describes the proposed fine-tuning for DeePC in detail. Section 4 validates the proposed approach through numerical simulations, and Section 5 concludes this work and discusses potential avenues for future research.

\textbf{Notation:} The initial and future values of \( e \) are denoted as \( e_{\mathrm{ini}} \) and \( e_{\mathrm{f}} \), respectively. The sets of past and future input trajectories are represented as \( U_{\mathrm{p}} \) and \( U_{\mathrm{f}} \), and the sets of past and future output trajectories are denoted as \( Y_{\mathrm{p}} \) and \( Y_{\mathrm{f}} \). The norm \( \Vert e \Vert_n \) is defined as \( \sqrt{e^Tne} \). The joint probability density function of random variables \( X_1 \) and \( X_2 \) is represented as \( pf(X_1, X_2) \), and the conditional probability density function of \( X_1 \) given \( X_2 \) is denoted as \( pf(X_1 | X_2) \). A random variable \( X \) following a Gaussian distribution with mean \( \mu \) and variance \( \theta \) is denoted as \( X \sim \mathscr{N}(\mu, \theta) \), whereas a uniformly distributed random variable with mean \( \mu \) and variance \( \theta \) is denoted by \( X \sim U(\mu, \theta) \). 

\section{Problem Formulation}

Consider a finite-dimensional, finite-horizon linear time-invariant system \( \mathcal{B} \) with input \( u(t) \) (\( u(t) \in \mathbb{R}^m, m \in \mathbb{N}^+ \)) and output \( y(t) \) (\( y(t) \in \mathbb{R}^p, p \in \mathbb{N}^+ \)), characterized by parameters of the number of inputs \( \mathbf{m}(\mathcal{B}) \), the lag \( \mathbf{l}(\mathcal{B}) \), and the order \( \mathbf{n}(\mathcal{B}) \).

\begin{equation}
	\mathcal{B}:\left[ \begin{array}{c}
		\sigma x(t)\\y(t)
	\end{array}\right]=\left[\begin{array}{cc}
	A&B\\C&D
	\end{array} \right] \left[ \begin{array}{c}
	x(t)\\u(t)
	\end{array}\right] 
	\label{system}
\end{equation}
where \( t \) is the time instant, with \( t = 1, 2, \ldots \), and \( A, B, C, D \) are matrices of appropriate dimensions, the shift operator \( \sigma \) is defined such that \( \sigma x(t) = x(t + 1) \).

\begin{figure}[!htb]
\centering%% For centre alignment of image.
\includegraphics[width=9cm,height=7cm]{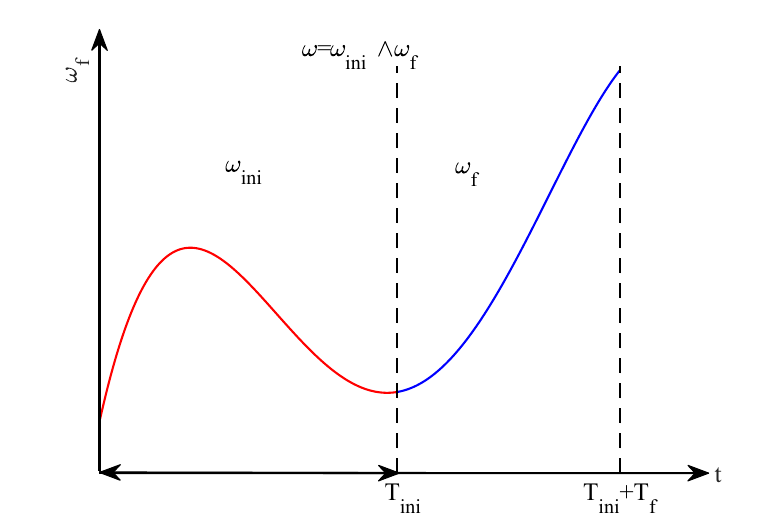}
\caption{Control from A DeePC Perspective. $\omega_{ini}$ and $\omega_f$ denote initial and future input-output behavior, respectively.}\label{ddc_control}
\end{figure}

We represent the systems \( \mathcal{B} \) satisfying equation \eqref{system} as \( \mathcal{L}^q \), where \( q = m + p \). The input-output behavior of the system at time \( T \) is represented by \( \omega(T) = [u(T) \quad y(T)]^T \in \mathbb{R}^{q \times 1} \). For the system described by \eqref{system}, DeePC views the system's behaviour as a collection of trajectories, each of which is defined as the set of all the input-output signals over a finite range \( [1, T] \):
\begin{equation}
	\omega|_{1:T}=\left( \omega(1),\ldots,\omega(T)\right)\in (\mathbb{R}^q)^T  
\end{equation}
Consider a linear time-invariant (LTI) system \( \mathcal{B} \in \mathcal{L}^q \) within the range \( [1, T] \)
\begin{equation*}
	\mathcal{B}|_T:=\left\lbrace \omega|_{1:T} \in \mathcal{B} \right\rbrace 
\end{equation*}
 %% Labels are used to cross-reference an item using \ref command.
DeePC methods formulate the control of a target system as identifying a future trajectory that satisfies specified requirements from all possible trajectories, as illustrated in Fig. \ref{ddc_control}. Note that if \( (T_1, \omega_1) \) and \( (T_2, \omega_2) \in \mathcal{B} \) with \( \omega_1(t_0) = \omega_2(t_0) \) at time \( t_0 \), then \( (T_1, \omega_1) \wedge (T_2, \omega_2) \in \mathcal{B} \), where the notation \( \wedge \) denotes the connection at time \( t_0 \).
\begin{equation*}
	\left(T_1, \omega _{1}\right) \wedge\left(T_2,\omega_{2}\right):=\left\{\begin{array}{ll}
		\left(T_1,\omega_{1}\right)(t) & \text { for } t<t_{0}, \\
		\left(T_2,\omega_{2}\right)(t) & \text { for } t \geqslant t_{0}
	\end{array}\right.
\end{equation*}

With a trajectory \( \omega_d \in (\mathbb{R}^q)^{T_d} \) that is long enough, we can construct a Hankel matrix with depth \( L \), where \( L \leq T_d \), as follows:
\begin{equation}
	\mathscr{H}_{L}\left(w_{\mathrm{d}}\right):=\left[\begin{array}{cccc}
		w_{\mathrm{d}}(1) & w_{\mathrm{d}}(2) & \cdots & w_{\mathrm{d}}\left(T_{\mathrm{d}}-L+1\right) \\
		w_{\mathrm{d}}(2) & w_{\mathrm{d}}(3) & \cdots & w_{\mathrm{d}}\left(T_{\mathrm{d}}-L+2\right) \\
		\vdots & \vdots & & \vdots \\
		w_{\mathrm{d}}(L) & w_{\mathrm{d}}(L+1) & \cdots & w_{\mathrm{d}}\left(T_{\mathrm{d}}\right)
	\end{array}\right] 
	\label{Hankel}
\end{equation}

The necessary and sufficient condition for utilizing the Hankel matrix to fully represent the system's behaviour is as follows:
\begin{equation}
	rank(\mathscr{H}_{L}\left(w_{\mathrm{d}}\right))=mL+n
	\label{fundamental_lemma}
\end{equation}
where $n$ is the minimal state dimension and the above equation is called the fundamental lemma. All trajectories of length \( T_d \) that satisfy this lemma can fully capture the characteristics of the system. In DeePC, \( \mathscr{H}_{L}(w_{\mathrm{d}}) \) is regarded as a non-parametric representation model drawn from the raw data from system \( \mathcal{B} \).

Consider a controllable system \( \mathcal{B} \) that satisfies equation \eqref{fundamental_lemma}. When \( L = T \), there exist non-unique vectors \( g \in \mathbb{R}^{T_d-T+1} \) such that
 \begin{equation}
 	\omega\in\mathcal{B}|_T\Leftrightarrow\mathscr{H}_{T}(\omega_d)g=\omega
 	\label{omega_represent}
 \end{equation}

The Hankel matrix includes the initial trajectory and the future trajectory as follows:
 \begin{equation*}
 	\mathscr{H}_{T_{ini}+T_f}(\omega_d)\thicksim\left[\frac{\begin{array}{c}
 			U_{\mathrm{p}} \\
 			U_{\mathrm{f}} \end{array}}{\begin{array}{c}
 			Y_{\mathrm{p}} \\
 			Y_{\mathrm{f}}
 	\end{array}} 
 	\right]=\left[\frac{\mathscr{H}_{T_{\text {ini }}+T_{f}}\left(u_{\mathrm{d}}\right)}{\mathscr{H}_{T_{\text {ini }}+T}\left(y_{\mathrm{d}}\right)}\right]
 \end{equation*} 
where $\thicksim$ represents the similarity of data arranged within a matrix.

However, due to noises in the raw data or the inherent system nonlinearity and complexity, the Hankel matrix \( \mathscr{H}_{L}(w_{\mathrm{d}}) \) is typically of full rank. As a result, any trajectory within the data-driven framework could be deemed a feasible solution, leading to potential ambiguity in trajectory selection. In practical implementations, obtaining noise-free data is challenging. Consequently, directly utilizing \( \mathscr{H}_{L}(w_{\mathrm{d}}) \) for optimal control may lead to suboptimal performance or the generation of infeasible control solutions.  

Though the raw data \( \omega \) collected from system \( \mathcal{B} \) contains uncertainties, DeePC employs regularization techniques to mitigate the impact of noisy measurements and accommodate mild nonlinearities. Coulson \textit{et al.} revealed that the DeePC framework is mathematically equivalent to model predictive control (MPC) in deterministic LTI systems. Given a reference trajectory \( \omega_r = (u_r, y_r) \in \mathbb{R}^{qT_f}, T_f > 0 \), DeePC incorporates \( l_1 \)-norm regularization through:
 \begin{equation}
 	\begin{array}{ll}
 		min&over\, g\;||\mathscr{H}_{T_{ini}+T_f}(\hat{\omega}_d^\star)g-\omega_{r,ini}||_P^2+\lambda_g\cdot ||g||_1\\
 		s.t.&\hat{\omega}_d^\star\in argmin\,over\,||\omega_d-\hat{\omega}_d||\\
 		&s.t.\quad rank \mathscr{H}_{T_{ini}+T_f}(\hat{\omega}_d)\leq m\cdot(T_{ini}+T_f)+n
 	\end{array}
 	\label{DeePC}
 \end{equation}
 where $\lambda_g\geq0$, $\omega_{r,ini}=(u_{ini},y_{ini},u_r,y_r)$. 
 
DeePC utilizes the most recent \( T_{\text{ini}} \) measurements to predict the system's future input-output behavior, so as to achieve optimal control based on the system formulation in \eqref{DeePC}. DeePC's primary objective is to determine the optimal control input that steers the system toward the desired output trajectory.  

In this work, we aim to establish a quantitative relationship between hyperparameters and the system's input-output behavior. We assume that the output of system \( \mathcal{B} \) is subject to additive Gaussian noise with mean \( \mu \) and variance \( \theta \), and the collected data satisfies the fundamental lemma condition in \eqref{fundamental_lemma}. The system dynamics are governed by \eqref{system}, and the control is implemented using the DeePC framework as formulated in \eqref{DeePC}. The following metric $M$ can be used to evaluate the impact of hyperparameter under stochastic disturbances: 

\begin{equation}
 	M=\sqrt{\frac{1}{n}\sum_{i=1}^{n}(y_{r,i}-\hat{y}_i)^2}
\end{equation} 
where \( n \) represents the number of samples, \( y_{r,i} \) denotes the reference trajectory, and \( \hat{y}_i = y_{c,i} + y_{d,i} \) is the actual control output, where \( y_{c,i} \) is the system's output under the input \( u_i \), and \( y_{d,i} \) is the environmental noise at time \( i \). The input energy is defined as the \( L_2 \)-norm of \( u \), denoted by \( \Vert u \Vert_{2,n} \).
 
\begin{figure}[!htb]
\centering%% For centre alignment of image.
\includegraphics[width=9cm,height=6cm]{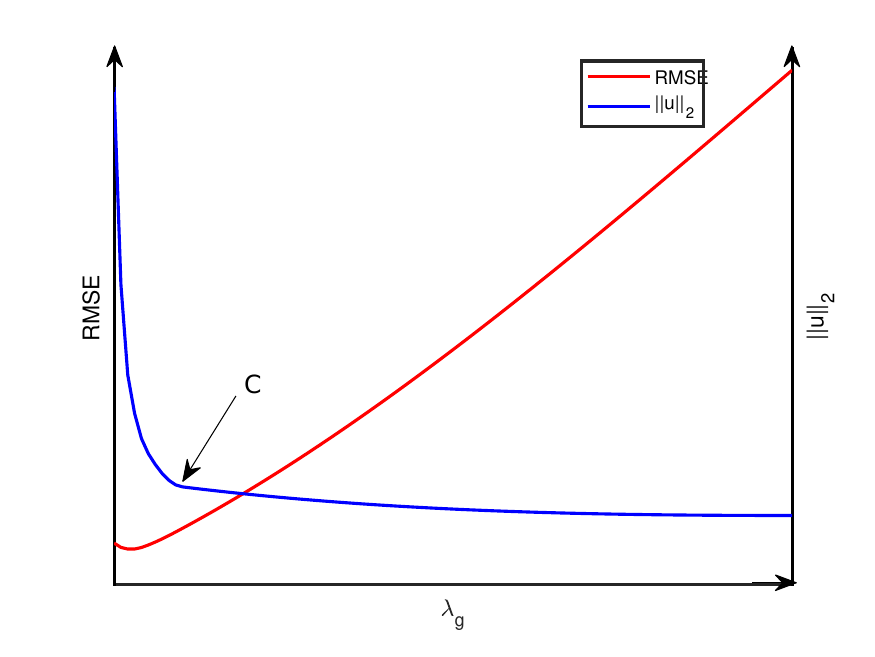}
\caption{Relationship between DeePC Control input Energy, M (measured by RMSE), and Hyperparameter $\lambda_g$.}\label{relationship}
\end{figure}

The interdependence among the control input norm \( \Vert u \Vert_{2,n} \), \( M \), and the hyperparameter \( \lambda_g \) is depicted in Fig. \ref{relationship}. Empirical analysis reveals that the increase in \( \lambda_g \) leads to rapid decline in \( \Vert u \Vert_{2,n} \), followed by a substantial reduction in the rate of decrease beyond a critical threshold, denoted as point \( C \). The behavior of \( M \) follows a non-monotonic trend, characterized by a slight initial decrease, after which it begins to rise. We also observe such a trend from Fig. \ref{relationship}, and the precise positioning of key inflection points may vary depending on system parameters and noise characteristics.

We would like to note that there lacks a concrete approach in tuning the hyperparameter \( \lambda_g \) in equation \eqref{DeePC} in response to uncertainties like noises in the system. In the following, we provide our solution for fine-tuning of DeePC in detail.
 
\section{Fine-tuning for DeePC}

\begin{figure}[!htb]
	\centerline{\includegraphics[width=9cm,height=6cm]{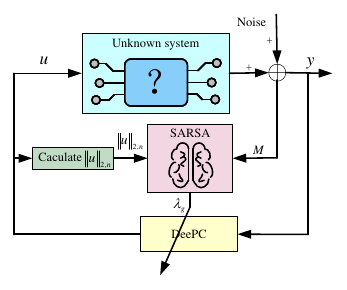}}
	\caption{Block diagram for the fine-tuning of DeePC.}
	\label{Q_DeePC}
\end{figure}

In this section, we introduce the fine-tuning of DeePC framework using reinforcement learning, which reformulates the problem of optimal hyperparameter selection in DeePC as a sequential decision-making process. This approach enables adaptive tuning of the regularization parameter in the control of systems with uncertainties and noises. The structural representation of the proposed fine-tuning of DeePC is depicted in Fig. \ref{Q_DeePC}. We detail our design in the following.

We define the following objective function when utilizing system I/O behavior
 \begin{equation}
 	J=|a\Vert u\Vert_{2,n}-bM|
 	\label{objective_function}
 \end{equation}
where the coefficients \( a \) and \( b \) are user-specified tuning parameters. In equation \eqref{objective_function}, those parameters regulate the trade-off between output error and control effort within the optimization objective. Specifically, when the control objective prioritizes minimizing output error, we can increase \( b \) and reduce \( a \), and vice versa when minimizing control effort is the primary concern. Analyzing Fig. \ref{relationship} suggests that an appropriate selection of \( a \) and \( b \) allows for the desired balance between these competing objectives, ensuring that equation \eqref{objective_function} can attain an optimal equilibrium.

We can reformulate the problem of determining the optimal parameter \( \lambda_g \) as follows:

\begin{equation}
 	\begin{array}{ll}
 		arg\;min&	J=|a\Vert u(\lambda_g)\Vert_{2,n}-bM(\lambda_g)|\\
 		s.t.&u\in \mathscr{U}\\
 		&y\in \mathscr{Y}\\
 		&\lambda_g\in \Lambda
 	\end{array}
	\label{problem1}
\end{equation}
where \( \mathscr{U} \) represents the set of feasible input values, \( \mathscr{Y} \) denotes the set of achievable output values of the system, and \( \Lambda \) is the set of system parameters. The constraints ensure that the solution to the problem remains feasible. In equation \eqref{problem1}, we assume that \( \Vert u(\lambda_g) \Vert_{2,n} \) and \( M(\lambda_g) \) are implicit functions of \( \lambda_g \). 

We would like to not that random additive noise in the collected data introduces significant challenges in modeling and analyzing the system's output behavior. The raw data affected by the noise is incorporated in the Hankel matrix formulation in \eqref{Hankel}. To address this issue, we leverage the probability distribution of the noise-corrupted signal as a decision variable for determining the optimal regularization parameter \( \lambda_{go} \). Specifically, we derive the conditional probability function \( pf(\lambda_g | \mathbf{x}) \) by discretizing the raw data and computing the joint probability distribution, where \( \mathbf{x} \) represents the system state.  

Once \( pf(\lambda_g | \mathbf{x}) \) is obtained, the optimal hyperparameter \( \lambda_{go} \) can be learned based on probabilistic inference. In this way, the computation of the objective function \( J \) at time step \( k \) relies on \( \lambda_{g,k} \) and the prior system state \( \mathbf{x}_{k-1} \), where \( \lambda_{g,k} \) denotes the regularization parameter at time \( k \) and \( \mathbf{x}_{k-1} \) represents the state at the preceding time step. To account for the influence of noise, a set of constraints is incorporated into the optimization problem defined in \eqref{problem1}. Upon the analysis above, we can reformulate \eqref{problem1} as a finite-dimensional stochastic optimization problem shown below, which can enable a robust and adaptive approach to hyperparameter adjustment
 
 \begin{subequations}
 	\begin{align}
 		arg\,min&J_k=|a\mathbf{x}_k(1)-b\mathbf{x}_k(2)|\\
 		s.t.&u\in \mathscr{U}\label{15b}\\
 		&y\in \mathscr{Y}\label{15c}\\
 		&\lambda_g\in \Lambda\label{15d}\\
 		&\mathbf{x}_{k} \sim pf_{\mathbf{x}, k}\left(\mathbf{x}_{k} \mid \mathbf{x}_{k-1}, \lambda_{g,k}=\hat{\lambda}_{g,k}\right)\forall k \in 1: T\label{15e}\\  
 		&\hat{\lambda}_{g,k} \sim pf_{k}\left(\hat{\lambda}_{g,k} \mid \mathbf{x}_{k-1}\right) \forall k \in 1: T \label{15f} \\
 		&\lambda_{g,k} \sim pf_{\lambda, k}\left(\lambda_{g,k} \mid \mathbf{x}_{k-1}\right) \forall k \in 1: T \label{15g}
 	\end{align}
 	\label{problem2}
 \end{subequations}
 
 where \( \mathbf{x}_k(1) \) represents the value of \( \Vert u \Vert_{2,n} \) at time step \( k \), \( \mathbf{x}_k(2) \) denotes the Root Mean Square Error (RMSE) \( M \) at time \( k \), and \( \lambda_{g,k} \) corresponds to the optimal hyperparameter selected at \( k \). The functions \( f_{\mathbf{u}, k}(\mathbf{u}_{k} \mid \mathbf{x}_{k-1}) \) and \( \mu_{k}(\hat{\mathbf{u}}_{k} \mid \mathbf{x}_{k-1}) \) belong to \( \mathcal{P} \), where \( \mathcal{P} \) denotes a probability space encompassing all possible probability distributions.  

Constraint \eqref{15e} enforces that the current system state \( x_k \) is probabilistically determined by the previous state \( x_{k-1} \), reflecting the stochastic nature of the system dynamics. Constraint \eqref{15g} specifies that at each time step \( k \), the optimal hyperparameter is inferred based on the probabilistic descriptions of the system's input and state dynamics.  

A key issue in solving the optimization problem is that, the probability space in equation \eqref{problem2} varies as the system or environmental conditions change. This variation complicates the optimization problem solving, since equation \eqref{problem2} requires explicit probability distributions. A valid solution necessitates solving the problem in \eqref{problem2} in an adaptive way, without prior knowledge about the probability space. It is also critical to ensure reasonably low computation in the optimization in facilitating real-time control.  

In addressing this challenge, we employ reinforcement learning to approximate the state-action of the system. Particularly, we adopt SARSA where we represent the system state-action by a Q-table~\cite{ma2022q}, and we update the table based on the interactions between the system and the environment.

\begin{figure}[!htb]
\centering%% For centre alignment of image.
\includegraphics[width=9cm,height=6cm]{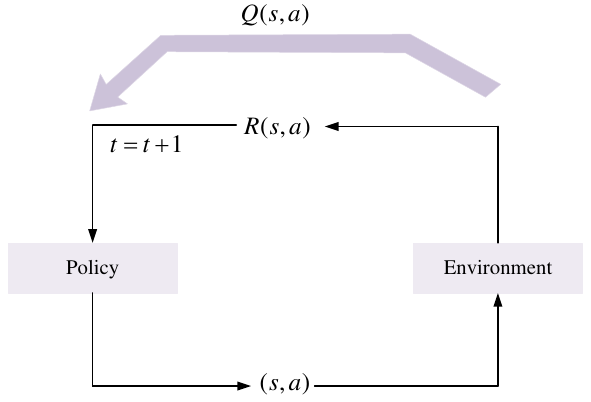}
\caption{SARSA architecture in DeePC tuning.}\label{SARSA}
\end{figure}

Fig. (\ref{SARSA}) shows the block diagram of the SARSA and illustrates an iterative learning process driven by continuous updates to the Q-table and policy. At each iteration, an action is selected based on the current state, executed, and subsequently evaluated as the system transitions to a new state. The set of feasible actions in the new state is then identified, and the Q-table is updated according to a predefined learning rule. Within this framework, when environmental changes occur, the reinforcement learning mechanism will update and retrain the Q-table to adaptively estimate the optimal DeePC hyperparameter, ensuring the control remains robust under varying system dynamics and uncertainties.

Under the reinforcement learning framework, equation \eqref{problem2} is a sequential decision-making problem. Specifically, in scenarios where the system's output is subject to Gaussian noise, the objective is to analyze the system state \( \mathbf{x} \), determine the optimal hyperparameter \( \lambda_{go} \), and maximize the cumulative reward \( R \) over time. We explain how we deal with the reinforcement learning in the following.

We collect data continuously during system operation, since the decisions made at the current time step not only affect the immediate system response but also have a cascading influence on future outputs. We aim to use reinforcement learning to dynamically adjust the DeePC hyperparameter in response to evolving system conditions
 \begin{equation}
 	\mathbf{x}_{k+1}=f_k(\mathbf{x}_k,\lambda_{g,k})
 	\label{zhuangtai}
 \end{equation}
where $f$ is a suitable function. Note that we do not need the explicit mathematical expression of $f$.
 
In equation \eqref{problem1}, the objective is to determine the optimal hyperparameter \( \lambda_{go} \) that minimizes the objective function \( J \) while simultaneously maximizing the cumulative reward \( R \). To this end, we formulate the reward function as the negation of \( J \) to ensure that the optimization framework aligns with reinforcement learning principles. The reward function is explicitly defined as follows:
 \begin{equation}
 	\begin{array}{ll}
 		R_k&=-|a\mathbf{x}_k(1)-b\mathbf{x}_k(2)|\\
 		&=g(\mathbf{x}_k,\lambda_{g,k})
 	\end{array}
 	\label{reward_function}
 \end{equation}
where $g(\cdot)$ is some function. To avoid symbol misuse, we will denote \( g(\mathbf{x}_k,\lambda_{g,k}) \) as \( R_k(\mathbf{x}_k,\lambda_{g,k}) \).
 
At this stage, the reward associated with each state-action pair can be explicitly determined. Consequently, equation \eqref{problem2} can be reformulated within the framework of reinforcement learning, where the objective is to derive an optimal policy that selects a sequence of hyperparameters \( \lambda_g \) to maximize the cumulative reward of the system over time.
 \begin{equation}
 	\begin{array}{ll}
 		max_{\lambda_{g,k}}&\mathbb{E}\left[\sum_{k=0}^{T}R_k\right]\\
 		s.t.&\mathbf{x}_{k+1}=f(\mathbf{x}_k,\lambda_{g,k})
 	\end{array}
 	\label{problem3}
 \end{equation}
 where $\mathbb{E}[\cdot]$ denotes expectation.

The objective is to determine an optimal policy, denoted as \( \pi_k \), that maximizes the expected reward. Note that the computation of \( \pi_k \) depends on the availability of historical state-action pairs. By taking \( \pi_k \) as the decision variable, the primary goal is to derive a function that optimally maps system states to the hyperparameter selection, expressed as \( \lambda_{g,k} = \pi_k(\mathbf{x}_k, \mathbf{x}_{k-1}, \ldots) \). Given an initial state \( x_0 \), the optimization problem can be reformulated as
 \begin{equation}
 	\begin{array}{ll}
 		max_{\pi_k}&\mathbb{E}\left[\sum_{k=0}^{T}R_k(\mathbf{x}_k,\lambda_{g,k})\right]\\
 		s.t.&\mathbf{x}_{k+1}=f(\mathbf{x}_k,\lambda_{g,k})\\
 		&\lambda_{g,k}=\pi_k(\mathbf{x}_k,\mathbf{x}_{k-1},\ldots)
 	\end{array}
 	\label{problem4}
 \end{equation}
 
Define the Q-function as
 \begin{equation}
 	\begin{array}{ll}
 		Q_{z\to c}(\mathbf{x},\lambda_g)&=\mathop{max}\limits_{\lambda_g}\;\mathbb{E}\left[\sum\limits_{k=z}^cR_k(\mathbf{x}_k,\lambda_{g,k})\right]\\
 		&\begin{array}{ll}
 			s.t.&\mathbf{x}_{k+1}=f(\mathbf{x}_k,\lambda_{g,k})
 		\end{array}
 		\label{Qfunction}
 	\end{array}
 \end{equation}
The Q-function represents the maximum value for equation \eqref{problem4} over the time from \( z \) to \( c \), given that the control input at time \( a \) is \( \lambda_g \) and the initial system state is \( \mathbf{x} \). After the optimal value of equation \eqref{problem4} is attained, the optimal policy is expressed as \( \pi(\mathbf{x}_0) = \arg \max_{\lambda_g} Q \to T(\mathbf{x}, \lambda_g) \), where the policy selects the action \( \lambda_g \) that maximizes the Q-value.  

In this formulation, only the current state of equation \eqref{problem4} is required to predict future state transitions, facilitating the application of dynamic programming for Q-function computation. Define the terminal Q-function as
 \begin{equation}
 	Q_{T\to T}(\mathbf{x},\lambda_g)=\mathbb{E}\left[R_T(\mathbf{x}_k,\lambda_{g,k})\right]
 \end{equation}
and calculate recursively at this stage
  \begin{equation}
 	Q_{k\to T}(\mathbf{x},\lambda_g)=\mathbb{E}\left[R_k+\mathop{max}_{\hat{\lambda}_g}Q_{k+1\to T}(f_k(\mathbf{x}_k,\lambda_{g,k}),\hat{\lambda}_g)\right]
  \label{Q_change}
 \end{equation}
and \eqref{Q_change} has the following deformation
 \begin{equation}
 	\begin{array}{l}	
 		Q_{k\to T}(\mathbf{x},\lambda_g)=\mathop{max}\limits_{\pi_{k+1},\ldots,\pi_T}\,\mathbb{E}\left[R_k+\sum\limits_{s=k+1}^{T}R_s(\mathbf{x}_s,\pi_s(\mathbf{x}_s))\right]\\
 		=\mathbb{E}\left[R_k+\mathop{max}\limits_{\pi_{k+1},\ldots,\pi_T}\mathbb{E}\left\{\sum\limits_{s=k+1}^{T}R_s(\mathbf{x}_s,\pi_s(\mathbf{x}_s))\right\} \right]\\
 		=\mathbb{E}\left[R_k+\mathop{max}\limits_{\pi_{k+1}}Q\left\{f(\mathbf{x},\lambda_g),\pi_{k+1}(f(\mathbf{x},\lambda_g))\right\}\right]\\
 		=\mathbb{E}\left[R_k+\mathop{max}\limits_{\hat{\lambda}_g}Q\left\{f(\mathbf{x},\lambda_g),\hat{\lambda}_g\right\}\right]
 	\end{array}
 	\label{Qchange}
 \end{equation}
 
In equation \eqref{Qchange}, substituting equation \eqref{Qfunction} yields the initial transformation, which establishes the relationship between the Q-function and the system's state-action dynamics. The final transformation provides a systematic approach for approximating the optimal hyperparameter at each iteration, facilitating its adaptive selection throughout the decision-making process.
 \begin{equation}
 	\lambda_{g,k}=\mathop{max}\limits_{\hat{\lambda}_g}Q_{k\to T}(\mathbf{x}_k,\lambda_g)
 	\label{outu}
 \end{equation}
Note that this estimate only depends on the current state. During the reinforcement learning, we update the Q-table in SARSA through

 \begin{equation}
 	\begin{array}{l}
		Q_{new}(\mathbf{x}_{k-1},\lambda_{g,k})\leftarrow Q_{old}(\mathbf{x}_{k-1},\lambda_{g,k})+\alpha(R_k+\\
	\gamma Q_{old}(\mathbf{x}_{k},\lambda_{g,k+1})-Q_{old}(\mathbf{x}_{k-1},\lambda_{g,k}))
	\label{train_function}
 	\end{array}
 \end{equation}
where \( \alpha \) represents the learning rate and \( \gamma \) denotes the discount factor. Equation \eqref{train_function} enables iterative updates without requiring explicit knowledge about the system's dynamic, thereby facilitating the identification of the optimal control action. In each iteration of the SARSA algorithm, the training dataset consists of a quintuple \( (\mathbf{x}_{k}, \lambda_{g,k}, R_k, \mathbf{x}_{k+1}, \lambda_{g,k+1}) \), where state-action transitions and rewards are updated sequentially according to the temporal progression of the system. At each update step, the next state-action pair is determined based on the current state-action pair, subsequently serving as the initial condition for the next iteration. 

To mitigate the risk of premature convergence to suboptimal solutions and prevent the algorithm from becoming trapped in local optima, an \(\epsilon\)-greedy policy is incorporated in the reinforcement learning. This policy enhances the learning efficiency of SARSA by maintaining a balance between exploration and exploitation, thereby improving the overall quality of the derived control strategy. The \(\epsilon\)-greedy policy is shown below

 \begin{equation}
 	\mu(\hat{\lambda}_{g,k}|\mathbf{x}_{k-1})=(1-\epsilon)\cdot pf_{\lambda_g}(\lambda_{g,k}|\mathbf{x}_{k-1})+\epsilon\cdot unif(U_k)
 	\label{updata_fun}
 \end{equation}
where $\epsilon$ is greedy constant, $\lambda_{g,k}\in \mathop{arg\,max}\limits_{\lambda_g}Q_{old}(\mathbf{x}_{k-1},\lambda_g)$,
 \begin{equation*}
 	pf_{\lambda_g}(\lambda_{g,k}|\mathbf{x}_{k-1})=\left\{\begin{array}{ll}
 		1&\lambda_{g,k}\in\Lambda\\0&\text{else}
 	\end{array}\right.
 \end{equation*}
and the function for the optimal parameter shifts to~\cite{andrew2018reinforcement}
 \begin{equation}
 	\lambda_{g,k}=\left\{\begin{array}{ll}
 		\mathop{max}\limits_{\hat{\lambda}_g}Q_{k\to T}(\mathbf{x}_k,\lambda_g)& \text{with probability} \;
 		1-\epsilon\\
 		\text{a random action in} \;\Lambda&\text{with probability}\; \epsilon
 	\end{array}\right.
 	\label{choose-lambda}
 \end{equation}
 where $(\cdot)$ as an optional action.

The above fine-tuning can be seamlessly integrated into control systems to enable real-time training of the reinforcement learning model. Note that this fine-tuning leverages operational data for continuous learning, and it operate independently without requiring prior knowledge about the system. To avoid online reinforcement learning the renders high computation burden, we structure the developed fine-tuning into two phases: offline training and online deployment. The online deployment uses the trained model to make real time decisions. We illustrate the overall procedures of the developed approach in \textbf{Algorithm}~\ref{alg:lambda}.
 \begin{algorithm}[!htb]
 	\caption{S-DeePC }\label{alg:lambda}
 	\begin{algorithmic}
 		\STATE 
 	
 		\STATE \textbf{Offline train stage}
 		\STATE \hspace{0.5cm} {\textsc{Input:}}$U_{\mathrm{p}},Y_{\mathrm{p}},\alpha,\gamma,a,b$
 		\STATE \hspace{0.5cm}{\textsc{Output:}}$\mathscr{H}_{L}$, Q-table
 		\STATE \hspace{1cm} Hankel matrix $\leftarrow$\eqref{Hankel}
 		\STATE \hspace{1cm} Initialize Q-table
 		\STATE \hspace{1cm} loop: At each time step $k$
 		\STATE \hspace{1.1cm}$\vert$ \hspace{0.4cm} $R_k\leftarrow$ \eqref{reward_function}
 		\STATE \hspace{1.1cm}$\vert$ \hspace{0.4cm} $\lambda_{g,k},\mathbf{x}_k\leftarrow$ Discretize $U_{\mathrm{p}},Y_{\mathrm{p}},\lambda_g$, using \eqref{zhuangtai}
 		\STATE \hspace{1.1cm}$\vert$ \hspace{0.4cm} 
 		Q-table $\leftarrow$\eqref{train_function}
 		\STATE \hspace{1cm} end
 		\STATE \textbf{Online stage}
 		\STATE \hspace{0.5cm} {\textsc{Input:}}$U_{\mathrm{p}},Y_{\mathrm{p}},U_{\mathrm{f}},Y_{\mathrm{f}},\alpha,\gamma,\epsilon,a,b,$ (Q-table)
 		\STATE \hspace{0.5cm}{\textsc{Output:}}$u_k$
 		\STATE \hspace{1cm} Hankel matrix $\leftarrow$\eqref{Hankel}
 		\STATE \hspace{1cm} Initialize Q-table (Load Q-table)
 		\STATE \hspace{1.1cm}$ \text{loop: At time step $k$ } $
 		\STATE \hspace{1.1cm}$\vert$ \hspace{0.4cm} $R_k\leftarrow$ \eqref{reward_function}
 		\STATE \hspace{1.1cm}$\vert$ \hspace{0.4cm} $\lambda_{g,k},\mathbf{x}_k\leftarrow$ Discretize current $\mathbf{x},\lambda_g$
 		\STATE \hspace{1.1cm}$\vert$ \hspace{0.4cm}$
 		 Q_{new}(\mathbf{x}_{k-1},\lambda_{g,k})\leftarrow$ \eqref{train_function}
 		 \STATE \hspace{1.1cm}$\vert$ \hspace{0.4cm}Choose $\epsilon$
 		\STATE \hspace{1.1cm}$\vert$ \hspace{0.4cm}$ \mu(\hat{\lambda}_{g,k}|\mathbf{x}_{k-1})\leftarrow$ \eqref{updata_fun}
 		\STATE \hspace{1.1cm}$\vert$ \hspace{0.4cm}$ \lambda_{g,k}\leftarrow$\eqref{choose-lambda}
 		\STATE \hspace{1.1cm}$\vert$ \hspace{0.4cm}\text{Take action} $\lambda_{g,k} $, \text{Run DeePC}
 		\STATE \hspace{1cm}end
 	\end{algorithmic}
 	\label{caculambda}
 \end{algorithm}

\section{Simulations and results}
This section conducts two sets of simulations. The first one entails a comprehensive sensitivity analysis on the algorithm's performance and check how the performance can be affected. The second one conducts a comparative assessment of the developed approach based on a hand-picked system. The last one demonstrate our approach in a real system of spring control. We implement all the numerical experiments in MATLAB 2024a, utilizing the OSQP optimization toolbox.

For the first two simulations, we consider a second-order system $\mathcal{B}$ with its system matrices as
\begin{equation}
	\begin{array}{cc}
		A=\left[\begin{array}{cc}
			0.49&4\\-0.066&1.5
		\end{array} \right]&B=\left[\begin{array}{cc}
			0.01&0\\0&0.01
		\end{array} \right]\\
		C=\left[\begin{array}{cc}
			1&0\\0&1
		\end{array} \right]&D=\left[\begin{array}{cc}
			0&0\\0&0
		\end{array} \right]
	\end{array}
	\label{controlled_system}
\end{equation}

During the data collection phase of the developed approach, the system input \( u_{1,2} \sim \mathscr{N}(0,10^{-3}) \), while the outputs are affected by additive, independently and identically distributed Gaussian noise \( \mathscr{N}(0,10^{-6}) \). We use 600 noisy input-output data points to construct the Hankel matrix, ensuring that the system input received sufficiently sustained excitation to satisfy the fundamental lemma~\eqref{fundamental_lemma}. The other parameters used were \( Q = 100 \), \( R = 1 \), \( T = 600 \), \( T_f = 12 \), \( a = 1 \), and \( b = 1 \). We compute the average value of the objective function \( J \) using 2200 individual \( J \) values. The comparative methods employed in this study are the baseline PeePC, Hanke \cite{lazar2022offset}, and DeePC-Hunt \cite{cummins2024deepc}.

\subsection{Analysis on the impact of the DeePC parameters}

\subsubsection{Data collection and analysis}

\begin{figure}[!htb]
	\centerline{\includegraphics[width=9cm,height=6cm]{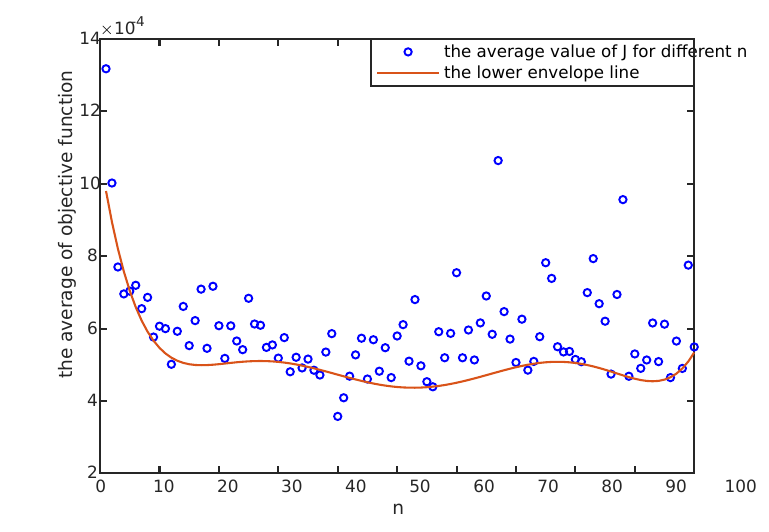}}
	\caption{Objective function under different values of parameter \( n \). }
	\label{n_yingxiang}
\end{figure}

We conducted a series of evaluations to determine the average value of the performance metric \( J \) across varying values of \( n \) within the range of 1 to 100, while keeping all other parameters constant. To facilitate trend visualization, a sixth-degree polynomial fitting was applied to construct the lower envelope, as depicted in Fig. (\ref{n_yingxiang}). In this figure, the blue scatter points represent the computed average values of \( J \) for different values of \( n \), while the red solid line denotes the lower envelope derived from polynomial fitting. Analysis of the scatter plot reveals a pronounced decline in \( J \) as \( n \) increases up to \( n = 13 \), beyond which the rate of decrease diminishes significantly, reaching a minimum at \( n = 40 \). Additionally, for \( n > 60 \), the distribution of \( J \) exhibits increased dispersion, which can be attributed to the system latency introduced by higher values of \( n \). This observation suggests that the performance metric \( J \) stabilizes, indicating diminishing returns for further increases in \( n \).

Considering the trade-off between computational delay and algorithmic performance, we identify \( n = 40 \) as the optimal choice. 

\subsubsection{Parameter selection}
\begin{figure}[!htb]
	\centerline{\includegraphics[width=9cm,height=6cm]{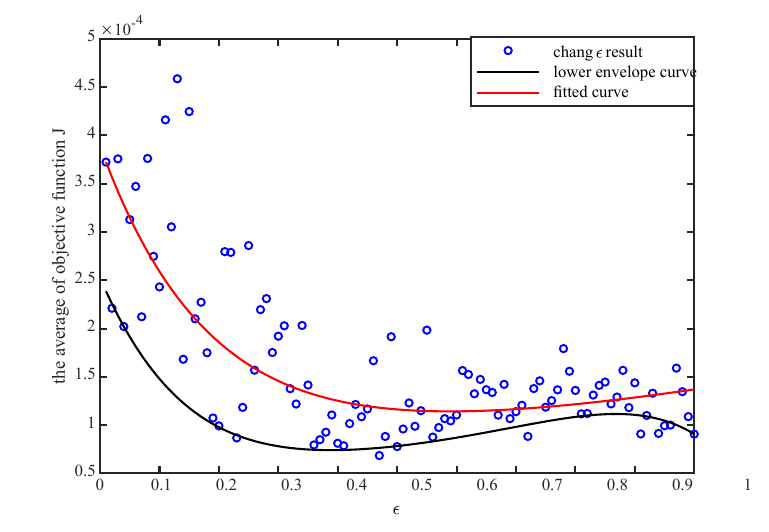}}
	\caption{Objective function under difference values of parameter $\epsilon$. }
	\label{eps_result}
\end{figure}

In the greedy exploration strategy, we investigated the influence of three parameters: the learning rate \( \alpha \), the discount factor \( \gamma \), and the exploration parameter \( \epsilon \). To isolate the effect of each parameter, we conducted a sensitivity analysis by varying one parameter while holding the others constant. For instance, with the learning rate fixed at \( \alpha = 0.5 \) and the discount factor at \( \gamma = 0.9 \), we examined the impact of \( \epsilon \) over the range \([0,1]\) with a step size of 0.01. The results were visualized using a scatter plot, which exhibited an exponential distribution pattern. The fitted model yielded a coefficient of determination \( R^2 = 0.94 \), and the lower envelope curve was obtained, as illustrated in Fig. \ref{eps_result}.  

In Fig. \ref{eps_result}, the blue scatter points represent the average objective function values corresponding to different values of \( \epsilon \), while the red solid line denotes the exponential fitting curve, and the black solid line depicts the lower envelope obtained via sixth-degree polynomial fitting. The fitted trend suggests that for \( \epsilon > 0.4 \), the rate of decrease in the objective function value becomes less pronounced. Furthermore, the lower envelope curve indicates that for \( \epsilon > 0.2 \), the decline in the objective function value slows considerably, and beyond \( \epsilon = 0.4 \), the average value increases. Further, the scatter distribution reveals that after \( \epsilon = 0.2 \), the lowest recorded objective function values demonstrate minimal variation. Based on these findings, we identified \( \epsilon = 0.35 \) as the optimal value.  Following the same methodology, the optimal values for the remaining parameters were determined to be \( \gamma = 0.86 \) and \( \alpha = 0.53 \).

\subsection{Fine-tuning performance}\label{simB}
We evaluate the developed algorithm's performance in comparison with the baseline DeePC, Hanke, and DeePC-Hunt. In all the cases, the parameter settings adhere to the previously determined optimal values to ensure consistency and comparability. We first check the convergence performance of those considered approaches, and then we evaluate those approaches under both Gaussian and uniformly distributed disturbance.

\begin{figure}[!htb]
	\centerline{\includegraphics[width=9cm,height=6cm]{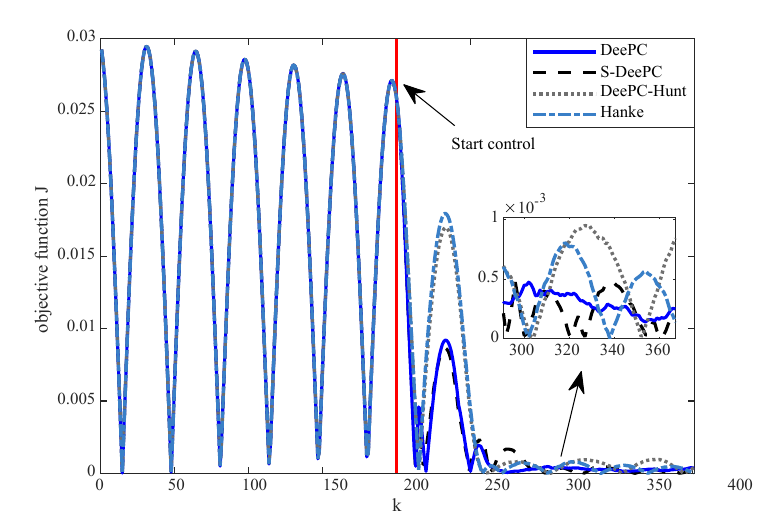}}
	\caption{Convergence comparison between Baseline DeePC, S-DeePC, DeePC-Hunt, and Hanke Algorithms. S-DeePC denotes the proposed approach. }
	\label{convergence}
\end{figure}

\subsubsection{Convergence speed}\label{Convergence_speed}
We analyze the convergence of the considered approaches by tracking the number of iterations required to achieve stability from its initial state, and we show the results in Fig. \ref{convergence}. We observe that, upon the initialization of the proposed approach, the rate of decline in the performance metric \( J \) accelerates, with the transient phase duration decreases from 16 steps to 12 steps before the oscillations. We observe the peak value of \( J \) at step 233, which is significantly lower than the corresponding value in the absence of control intervention. Stability was attained at step 289, with \( J \) stabilizing within the range of \( 0 \) to \( 4.5 \times 10^{-4} \). The algorithm required a total of 89 training steps to achieve convergence.  

In comparison, the basedline DeePC was implemented with a fixed hyperparameter setting of \( \lambda_g = 0.03 \). While it demonstrated a similar initial transient response to the proposed S-DeePC, its oscillation amplitudes are more pronounced, although it exhibited a slightly faster convergence rate. The baseline DeePC reaches steady-state at step 264, attributed to the predetermined hyperparameter \( \lambda_g \), whereas S-DeePC adaptively optimized its hyperparameter during training.  

The Hanke and DeePC-Hunt methods exhibit slightly higher peak values of \( J \) than the baseline DeePC. However, both algorithms achieve a stable state at step 256. As illustrated in the inset of Fig. \ref{convergence}, though the Hanke and DeePC-Hunt methods converged earlier, their steady-state values of \( J \) remain higher than those achieved by both the baseline DeePC and the proposed S-DeePC approach, indicating suboptimal steady-state performance.

 \begin{figure*}[!htb]
 	\centerline{\includegraphics[width=20cm,height=6cm]{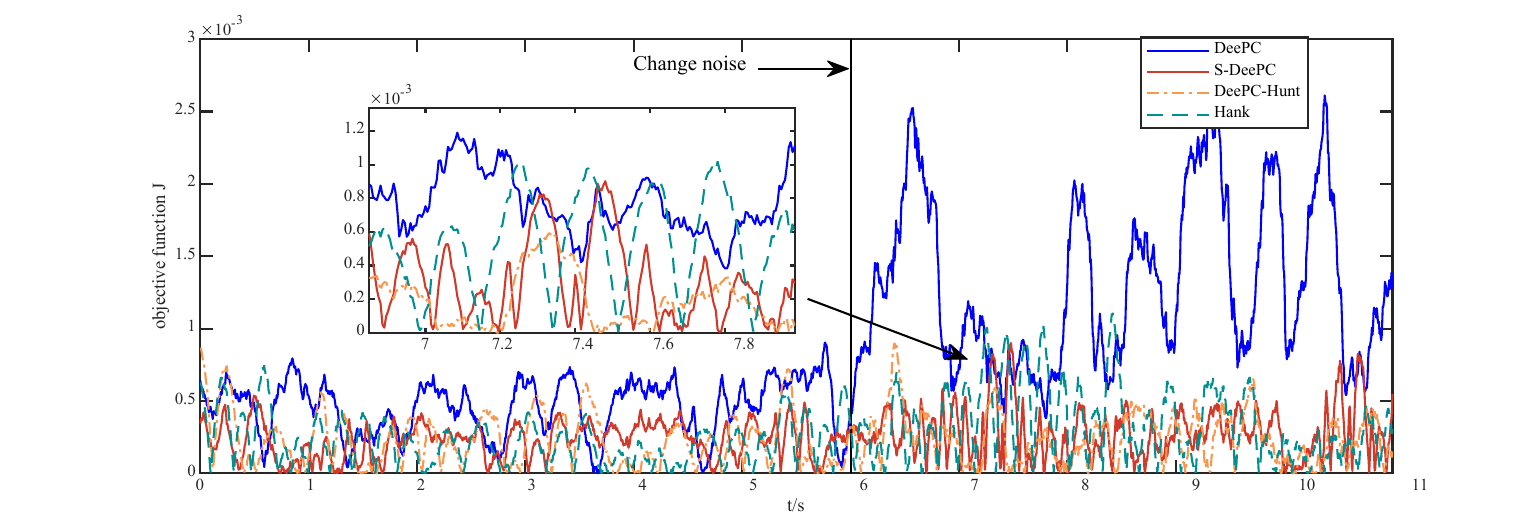}}
 	\caption{Control performance comparison between baseline DeePC, S-DeePC, DeePC-Hunt, and Hanke under normal and doubled Gaussian noise conditions.}
 	\label{evaluate}
 \end{figure*}

\subsubsection{Control of system under Gaussian distributed disturbance}\label{Gaussian_distribution}

In this experiment, a pre-trained Q-table is utilized to facilitate adaptive hyperparameter selection. The training dataset is generated under fixed system parameters, with \( \lambda_g \) incremented in discrete steps of 0.006 over the range \([0.006, 0.606]\), yielding a total of 1000 data points per \( \lambda_g \) value. At the onset of training, the Q-table is initialized as a zero matrix of dimensions \( 101 \times 101 \times 101 \), and the learning process follows the update rule defined in equation \eqref{train_function}.  

The system is initially subjected to Gaussian noise with a distribution of \( \mathscr{N}(0,10^{-6}) \). At \( t = 6 \) seconds, the noise energy is doubled to evaluate the robustness of the algorithms. The corresponding output trajectories for different control strategies are illustrated in Fig. (\ref{evaluate}).

During the first six seconds, we observe that both the baseline DeePC and S-DeePC maintain relatively stable performance, with S-DeePC exhibiting improved robustness compared to the baseline. However, upon the increase in noise intensity at \( t = 6 \) seconds, the output metric \( J \) of the baseline DeePC experiences a pronounced surge, followed by significant fluctuations, with most values exceeding the maximum \( J \) observed during the initial six-second period. In contrast, while the outputs of the S-DeePC, Hanke, and DeePC-Hunt also exhibit an increase, S-DeePC demonstrates a higher frequency of adaptive adjustments, mitigating extreme deviations and exhibiting only a marginal increase in peak values.  Further analysis of the locally magnified time window from \( t = 7.2 \) to \( t = 8 \) seconds reveals that under the intensified noise conditions, the performance of Hanke deteriorates, whereas S-DeePC and DeePC-Hunt maintain comparable performance levels. Beyond \( t = 8 \) seconds, the outputs of all the methods except the baseline DeePC progressively converge, attaining a steady-state within approximately 10 seconds.

\subsubsection{Uniformly Distributed Disturbance}
\begin{figure*}[!htb]
 	\centerline{\includegraphics[width=20cm,height=6cm]{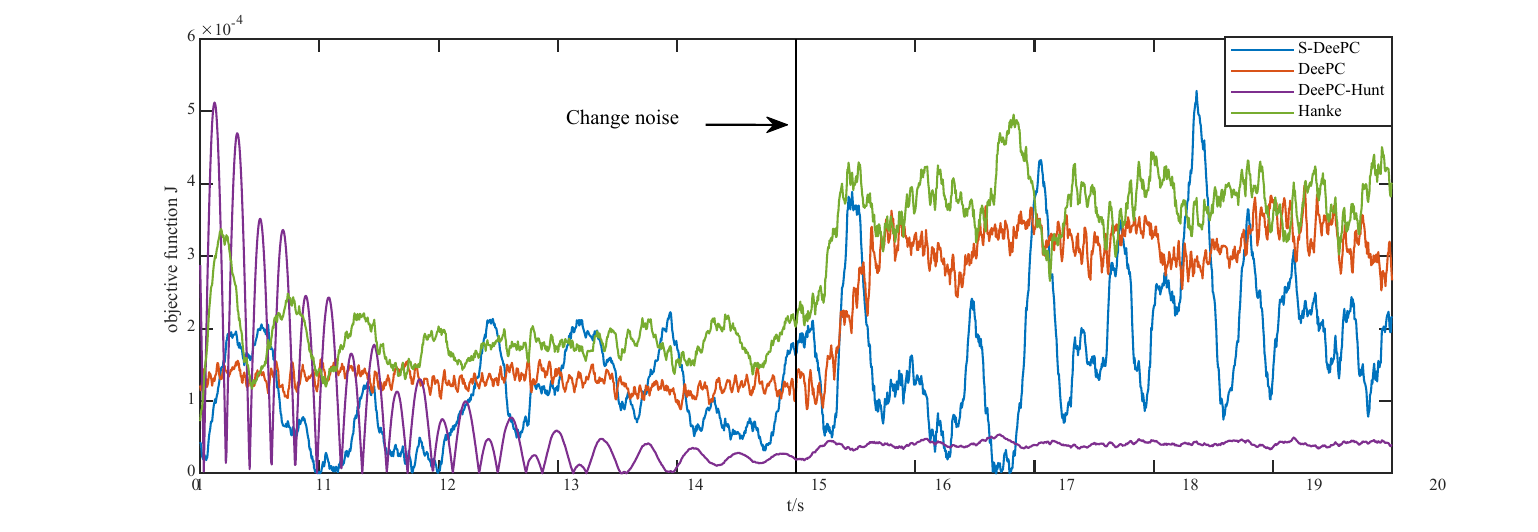}}
 	\caption{Control performance comparison between baseline DeePC, S-DeePC, DeePC-Hunt, and Hanke under normal and doubled uniform  noise conditions.}
 	\label{UniformlyDistributed}
 \end{figure*}

This simulation assumes that the disturbance follows \( U(0, 1 \times 10^{-5}) \). All algorithmic parameters are maintained consistent with those in Section~\ref{Gaussian_distribution} to ensure comparability. We collect training data under uniformly distributed disturbances and use the data for the assessment. we show the results in Fig. \ref{UniformlyDistributed}. At \( t = 15s \), the disturbance distribution transitions to \( U(0, 2 \times 10^{-5}) \), introducing an increase in noise intensity. We illustrate how this change affects the considered approaches in the following.  

As demonstrated in Fig. \ref{UniformlyDistributed}, S-DeePC generally outperforms the baseline DeePC, both before and after the noise variation. However, at certain peak values, its performance slightly lags behind the baseline DeePC, whereas the Hanke method exhibits a response pattern more closely aligned with the baseline DeePC. The Baseline DeePC maintains relatively stable output behavior, while S-DeePC demonstrates higher fluctuation amplitudes, reflecting its adaptive nature.  

A notable observation is that DeePC-Hunt achieves stability after approximately 14 seconds. Although its objective function value \( J \) increases momentarily following the change in noise, it exhibits the smallest fluctuation range among all tested algorithms. Additionally, the objective function value \( J \) for S-DeePC exhibits a gradual decreasing trend over time, indicating its ability to dynamically adjust and improve performance in response to changing noise conditions.
 
\subsection{Spring control study}
\begin{figure}[!htb]
	\centerline{\includegraphics[width=8cm,height=4cm]{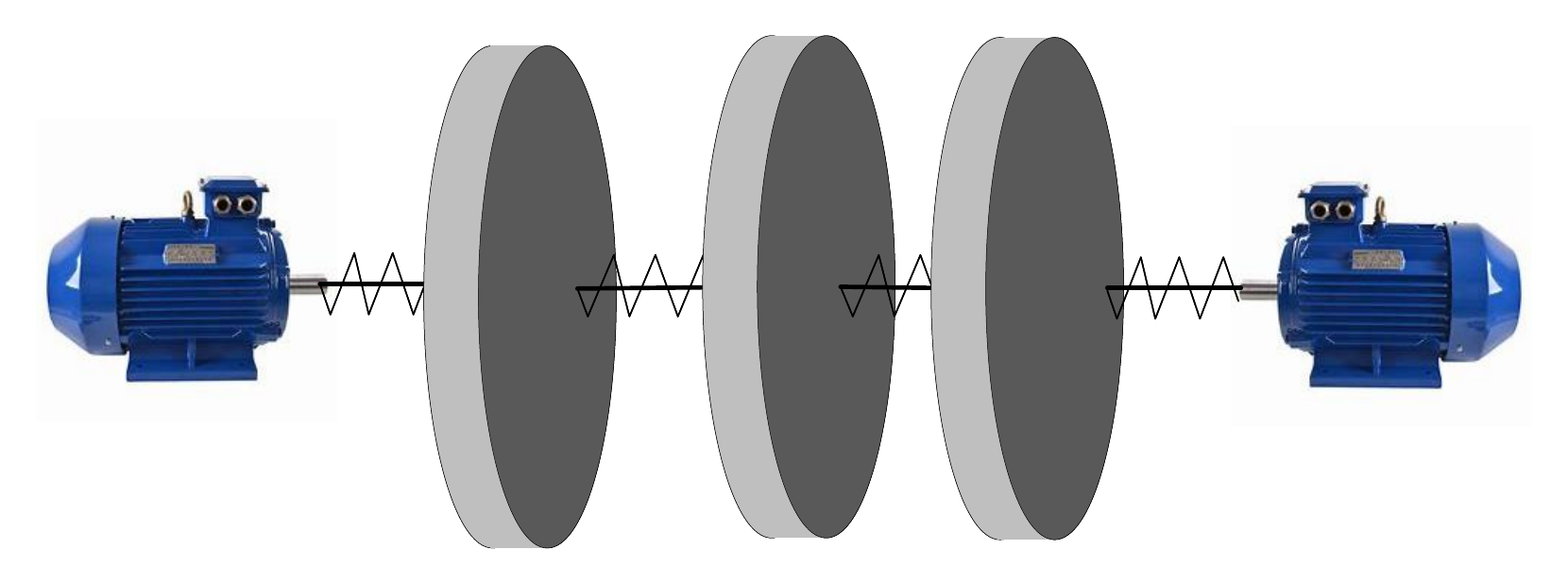}}
	\caption{The triple-mass-spring system with two motors and three discs.}
	\label{Triple-mass-spring}
\end{figure}

In this section, we demonstrate the proposed method using a Triple-mass-spring system~\cite{fiedler2021relationship} as shown in Fig. \ref{Triple-mass-spring}. This system consists of two stepper motors, springs, and three discs. Specifically, two motors are connected directly to the outermost discs, while the remaining discs are interconnected by springs. The angular positions of the discs serve as the system outputs. Thus, we have the parameters for this system that \( m = 2 \), \( P = 3 \), and the system lag is \( l = 2 \). The input constraints are set to the range \([-0.7,\,0.7]\). Data collection is conducted independently using the same parameter settings as in the preceding experiment. Gaussian noise with a variance of \(2.5\times10^{-9}\) is added to the system output to emulate sensor measurement noise.

For the proposed S-DeePC method, the Q-table is trained offline by sampling \(\lambda_g\) in increments of 0.01 over the interval \((0,1]\), and the Q-table is constructed once data collection has concluded. Likewise, Hanke’s algorithm is executed offline a single time, and the resulting \(\lambda_g\) is subsequently employed in the evaluations.

\begin{figure}[!htb]
 	\centerline{\includegraphics[width=9cm,height=6cm]{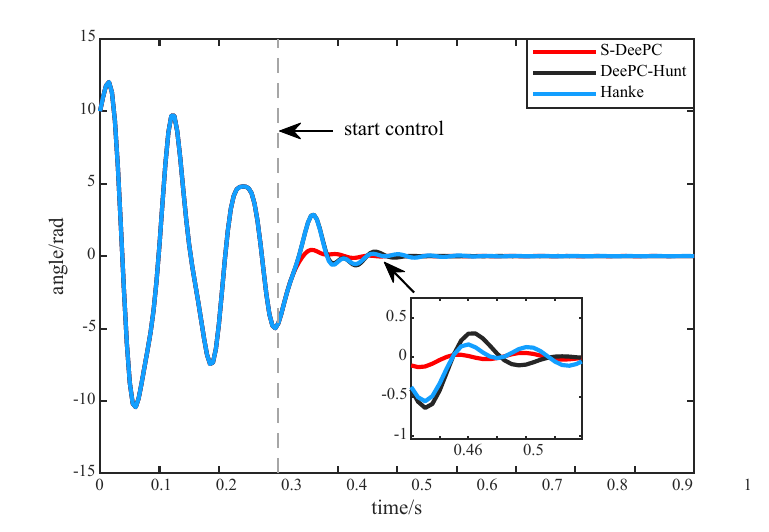}}
 	\caption{Control for the middle disc of the triple-mass-spring system. The control starts at 0.3 seconds.}
 	\label{UniformlyDistributed_spring}
 \end{figure}

In the initial phase of this simulation, all three discs are initialized at an angular displacement of 10 rad. At \(t = 0.3\), control inputs from both stepper motors were applied with the objective of returning the system to its equilibrium position. We focus on the central disc—chosen because it exhibits the largest dynamic lag. We plot its angular response in Fig. \ref{UniformlyDistributed_spring}, where we compare three approached: the proposed S-DeePC (red), DeePC-Hunt (black), and the Hanke method (blue). Prior to the control, the system undergoes passive damping due to the spring. After \(t = 0.3\) s, both motors start the control and drive the system, but the central disc responds with a slight delay. 

The proposed S-DeePC achieves stabilization with approximately 0.38 s, with peak-to-peak oscillations contained within \(\pm 0.15\) rad. In comparison, both DeePC-Hunt and Hanke require roughly 0.50 s to stabilize. As depicted in the inset graph, DeePC-Hunt outperforms Hanke method marginally during the interval \(0.30\)–\(0.48\) s. Subsequently, the performance relationship reverses, with DeePC-Hunt algorithm outperforming the Hanke method in the control.

\subsection{Summary of the simulation results}
 
The comparative results indicate that under noise conditions, Hanke and DeePC-Hunt methods achieve rapid convergence but exhibit relatively worse performance in the control. When the noise energy fluctuates, all the three adaptive approaches, \textit{i.e.}, Hanke, DeePC-Hunt, and the proposed S-DeePC eventually regain stability after a transient adaptation period. DeePC-Hunt method requires longer time to converge but demonstrates superior robustness to noises. Hanke method needs to recompute the optimal regularization parameter \( \lambda_g \) whenever the noise environment changes, which results in a computational delay of approximately 2–3 seconds. Compared with the two adaptive approaches above, the proposed S-DeePC algorithm exhibits faster convergence while maintaining competitive performance.

\section{Conclusion}\label{section5}

This work proposed a fine-tuning for DeePC for the control of uncertain systems using reinforcement learning. We analyzed the impact of DeePC hyperparameter on system input–output behavior, and reformulate the fine-tuning of DeePC hyperparameter as a sequential decision-making problem. We solved this problem through reinforcement learning, where we optimize DeePC hyperparameters to fit the varying system dynamics. Simulation results demonstrate that our method can effectively identify optimal DeePC hyperparameter in dynamic settings, and maintain stable control performance when disturbance conditions shift. We envision future work to test our fine-tuning in large-scale systems in terms of fine-tuning efficiency and real time performance. We also suggest that autonomous mechanism can be adopted to detect changes disturbances and system dynamics as a prognosis approach for the fine-tuning.

\bibliographystyle{unsrt}
\bibliography{myrefs}

\end{document}